\documentclass[11pt,psamsfonts,leqno,twoside]{amsart}
\usepackage{amsmath, amssymb}
\usepackage{amsfonts, setspace}
\usepackage{amsthm}

\allowdisplaybreaks
\newtheorem{theorem}{Theorem}[section]
\newtheorem{lemma}[theorem]{Lemma}
\newtheorem{proposition}[theorem]{Proposition}
\newtheorem{corollary}[theorem]{Corollary}

\newtheorem*{theorem*}{Theorem}
\theoremstyle{remark}
\newtheorem{remark}[theorem]{Remark}
\newtheorem{definition}[theorem]{Definition}

\numberwithin{equation}{section}

\def\Tr{{\rm Tr}}
\def\span{{\rm span}}
\def\co{{\rm co}}
\def\II1{II$_1$}
\def\unitary{\mathcal{U}}
\def\ignore#1{\relax}

\begin{document}

\title[On the normalizing algebra of a MASA in  a \II1 Factor ]
{On the normalizing algebra of a MASA in  a \II1 Factor }

\author{Ionut Chifan}
\address{ Department of Mathematics\\ University of Iowa\\ Iowa
City, Iowa 52242}
\email{ ichifan@math.uiowa.edu}

\maketitle

\begin{abstract}  Let $A$ be a maximal abelian subalgebra (MASA) in a \II1  factor $M$. Sorin Popa introduced an analytic condition
that can be used to identify the normalizing algebra of $A$ in  $M$ and which we call \emph{the relative weak asymptotic homomorphism property}. In this paper we show this property is always satisfied by the normalizing algebra of $A$ in $M$ and as a consequence we obtain that $\overline{\bigotimes}_{i\in I}(\mathcal{N}_{M_{i}}(A_{i})^{\prime\prime})= (\mathcal{N}_{\overline{\bigotimes}_{i\in I}M_{i}}(\overline{\otimes}%
_{i\in I} A_{i}))^{\prime\prime}$ .

\end{abstract}

\section{Introduction}

Throughout this paper $M$ will denote a fixed type \II1 factor and $A$ will denote a maximal abelian subalgebra (MASA) of $M$. By $\mathcal{U}(M)$ we denote the unitary group of $M$ and by $\mathcal{P}(M)$ the set of projections of $M$. Also, $\mathcal{N}_{M}(A)$ will denote the group of unitaries $u$ in $M$ such that $uAu^* = A$, i.e., $\mathcal{N}_{M}(A)$ is the \emph{normalizer} of $A$ in $M$. The von Neumann algebra $\mathcal{N}_{M}(A)^{\prime \prime}$ generated by $\mathcal{N}_{M}(A)$ will be called \emph{the normalizing algebra of $A$ in $M$}. We thus have $$
A \subseteq \mathcal{N}_{M}(A)^{\prime \prime} \subseteq M.$$  

A consistent study of algebra $\mathcal{N}_{M}(A)^{\prime \prime}$ started with the work of Jacques Dixmier. In  ~\cite{Di} the author distinguished three cases of particular interest, acording to the "size" of $\mathcal{N}_{M}(A)^{\prime \prime}$ in $M$:
\begin{enumerate}
\item \emph{Regular(or Cartan)} if $M= \mathcal{N}_{M}(A)^{\prime \prime}$.
\item \emph{Semi-regular} if $ \mathcal{N}_{M}(A)^{\prime \prime}$ is a \II1 factor.
\item \emph{Singular} if $A= \mathcal{N}_{M}(A)^{\prime \prime}$.
\end{enumerate}

Our primary objective in this paper is to show in Theorem \ref{2.1} that the triple $$
A \subseteq \mathcal{N}_{M}(A)^{\prime \prime} \subseteq M$$ has the following property:

\begin{definition}\label{1.1}
A triple of von Neumann algebras $B \subseteq N \subseteq M$, where $M$ is a \II1 factor, is said to have \emph{the relative weak asymptotic homomorphism property (relative WAHP)} when \newline for all  $x_{1},x_{2},x_{3},...x_{n} \in M$ and for every  $ \varepsilon
>0$ there exists  $u \in \unitary(B)$ such that 
$$\left\Vert
E_{B}(x_{i}ux_{j})-E_{B}(E_{N}(x_{i})uE_{N}(x_{j}))\right\Vert
_{2}<\varepsilon $$ for $i,j=1,\dots,n$. By $E_B$ and $E_N$ we denote the conditional expectations of $M$ onto $B$ and $N$, and by $\Vert x \Vert_{2}$ we denote the norm of $x$ taken in $L^2(M)$ - the Hilbert space of the GNS construction based on the faithful normal trace of the factor $M$.
\end{definition} 
In ~\cite{SiSm} the authors introduced the notion of \emph{weak asymptotic homomorphism property (WAHP)} referring to a MASA, $A\subset M$. We can immediately see this is equivalent to saying that the triple $A\subseteq A\subseteq M$ has the \emph{relative WAHP}, so our Definition \ref{1.1} can be viewed as a generalization of their concept.

Remarkably, in ~\cite{SiSmWhWi} it was shown that every singular MASA $A\subset M$ in a \II1 factor satisfies the WAHP, so our Theorem \ref{2.1} is also a generalization of their result. 

We would like to point out that Theorem \ref{2.1} represents a suitable tool to compute the normalizing algebra of a MASA in certain situations of tensor products of \II1 factors: 
\begin{corollary}

Let $I$ be a countable set, $\{M_i\}_{i\in I}$ a collection of $II_1$ factors and  $\{A_i\}_{i\in I}$ a collection of abelian diffuse von Neumann algebras, such that for every $i\in I$ we have $A_i\subset M_i$ is a MASA. Then,
 
 $$\overline{\bigotimes}_{i\in I}(\mathcal{N}_{M_{i}}(A_{i})^{\prime\prime})= (\mathcal{N}_{\overline{\bigotimes}_{i\in I}M_{i}}(\overline{\otimes}%
_{i\in I} A_{i}))^{\prime\prime}.$$
\end{corollary}

The same method can be used to estimate the normalizing algebra of certain subalgebras in cross-products (see Theorem \ref{3.8}).\newline

In order to study singular MASAs, Sinclair and Smith
~\cite{SiSm}
isolated a concept which they called \emph{strong singularity}:
\begin{definition}A MASA $A$ in a type \II1 factor $M$ is called \emph{strongly singular}  
   if  for every $u \in \unitary(M)$, 
$$\left\Vert E_{u{A}u^{\ast
}}-E_{A}\right\Vert _{\infty ,2}\geq _{_{{}}} \left\Vert
u-E_{A}(u)\right\Vert _{2}.$$
\end{definition}

By a very beautiful argument, Lemma 2.1 from ~\cite{RoSiSm} shows that the WAHP of $A\subset M$ a MASA implies the strongly singularity of $A$ in $M$. Since we have in hand Theorem \ref{2.1}, the same argument shows the following:  

\begin{theorem} 

 (A generalization of Sinclair-Smith inequality) Let $M$ be a  $II_{1}$ factor with $A\subset M$ a MASA. Then,

$$||E_{uAu^{\ast }}-E_{A}||_{\infty ,2}\geq ||u-E_{\mathcal{N}_{M}(A)^{\prime \prime}}(u)||_{2}, \:{ for}\: {all }\: u\in \mathcal{U}(M).$$

\end{theorem}

Here we would like to mention that even if the previous theorem is a generalization of strong singularity concept for a MASA in a \II1 factor, a more general version (with absolute constant 1) for Theorem 6.2 in ~\cite{PoSiSm} remains open.\newline

Even though the computations of the normalizing algebra of an arbitrary fixed subalgebra in a $II_{1}$ factor have proved to be a very difficult problem, in the last years we witnessed a constant and successful effort in this direction.
Without making a formal definition, Popa (~\cite{Po3,Po4,Po8}) and independently Robertson-Sinclair-Smith (~\cite{SiSm,RoSiSm}) verified WAHP for certain inclusions of von Neumann algebras $P\subset M$ and obtained containment of the normalizing algebra of $P$ in $M$ in various contexts. The case of free products and the case of weakly mixing actions of groups on von Neumann algebras are only few examples in this sense.

Furthermore, the ground breaking technology that Popa developed to control the normalizing algebra (and relative commutants in particular) works in a much more general setting ( ~\cite{Po3,Po4}), involving intertwining elements between two distinct subalgebras rather than normalizing elements of a common subalgebra. More precisely, Theorem 2.1 and Corollary 2.3 in ~\cite{Po4} give a complete description of the existence of intertwining elements between two distinct subalgebras in a fixed factor with discrete decomposition. This was called \emph{intertwining techniques} and was one of the major ingredients that has led to many striking results in von Neumann Algebras theory/Noncommutative Ergodic theory(~\cite {Po3,Po4,Po8,IPePo}).

However, in the situation $A\subset M$ is a MASA our Theorem \ref{2.1} shows that relative WAHP intrinsically characterizes the triple $A \subseteq \mathcal{N}_{M}(A)^{\prime\prime} \subseteq M$. 
 Our proof that the triple $A\subseteq  \mathcal{N}%
_{M}(A)^{\prime\prime} \subseteq M$ satisfies the relative WAHP is based on a deep idea of S. Popa ~\cite{Po3,Po4}, which is
to build normalizing elements by looking at the relative commutant between
the MASA and the basic construction ~\cite{Ch,Jo,Po1} for the inclusion $A \subseteq M$. This connection has been made before with great benefit, but we refine it.
Namely, we analyze the relationship between finite trace projections
in the the basic construction $\langle M, e_A \rangle$ 
 and the Jones projection $e_{N}$ of the normalizing algebra  $N$. It turns out that they satisfy an
interesting geometric relation that is revealed in the proof of Proposition 2.6. This is the key observation of our proof.

\bigskip

This paper is organized in two sections. In the first section we prove Theorem \ref{2.1}, which is the main result of the paper. In the second section we 
present some immediate applications of this theorem described above (Proposition \ref{3.3}, Proposition \ref{3.5}, Corollary  \ref{3.6},
Theorem \ref{3.8}).

\medskip

\emph{Acknowledgments:}

\bigskip

I would like to express my gratitude to my advisers, Paul Muhly and Florin Radulescu,
for constant guidance and illuminating discussons about this
material. I would also want to thank Fred Goodman for his enthusiastic
encouragement during the last years. Last but not the least I want to thank the referee and Prof. Roger Smith for the generous comments regarding the original version of this paper.

\section{An analytic characterization of normalizing algebra of a MASA in a $II_1$ factor}

First we recall some elementary properties of the basic construction for an
inclusion of two von Neumann algebras. These were thoroughly developed in
~\cite{Jo,Ch,Po1,PiPo}. Then in Theorem \ref{2.1} (which is the core result of this section) we provide an alternative description of
the absence of normalizing elements for a MASA $A$ in a \II1 factor $M$ which we call \emph{the relative weak asyptotic homomorhpism property(WAHP)}.

 Consider a triple $A\subseteq N\subseteq M$ , where $M$ is $II_{1}$
factor with the normalized trace $\tau $, $A$ a MASA, and $N$ is an intermediate von Neumann subalgebra.
The trace $\tau $ induces an inner product on $M$:
         $$\langle x,y\rangle=\tau(y^{*}x)\quad\text{for all}\quad x,y\in M $$

We denote by $L^{2}(M)$ the completion of $M$ with respect to the norm
$\|x\|_{2}=(\tau(x^{*}x))^{1/2}$,
and when an element $x \in M$ is seen as a vector in $L^{2}(M)$ it will be denoted by $\hat{x}$. 
Next, consider the $\tau$-preserving conditional expectations $E_{A}$ and $E_{N}$ onto $A$ and $N$, respectively.
Each such conditional expectation can also be viewed as a projection in $B(L^{2}(M))$. Thus we can define:\newline
$$e_{A}(\hat{x})=\widehat{E_A(x)}$$ for all $x \in M$. \newline Also there is an anticonjugation 
$\mathcal{J}:L^{2}(M)\rightarrow L^{2}(M)$, defined by
$$\mathcal{J}(\hat{x})=\widehat{x^*}$$ for $x \in M$.\newline Perform the basic construction with respect to $A$:
$$ A\subset M\subset \langle M,e_{A}\rangle \subset B(L^{2}(M))$$
 where  $\langle M,e_{A}\rangle$ is the von Neumann algebra generated by $M$ and $e_{A}$ in $B(L^{2}(M))$, i.e.   $\langle M,e_{A}\rangle=\{M,e_{A}\}^{''}$.

In the next proposition we recall several important properties of this construction that are of essential use to our next proofs. For the reader convenience we also include references that provide detailed proofs and a complete history of these facts.

\begin{proposition}\label{2.8}
Suppose $A\subseteq N\subseteq M$ is a triple of algebras as before. Then the following properties hold true:
\begin{enumerate}
\item [a.~\cite{Jo}] We have $\langle M,e_{A}\rangle=\mathcal{J}A^{\prime }\mathcal{J}$, which in particular says that $\langle M,e_{A}\rangle$ is a von Neumann algebra of
type $I_{\infty }$.

\item [b.~\cite{Jo}] For all $ x\in M$ we have $e_{A}xe_{A}=E_{A}(x)e_{A}$.  Moreover,  $\overline{\span}%
^{w}(Me_{A}M)=\langle M,e_{A}\rangle$,  and if $N=\mathcal{N}_{M}(A)$, then $e_{N}\in A\vee \mathcal{J}A\mathcal{J}$.
( by $A\vee \mathcal{J}A\mathcal{J}$ we denote the von Neumann algebra generated by $A$ and $\mathcal{J}A\mathcal{J}$
and by $\overline{\span}%
^{w}(Me_{A}M)$ the $w$-closure of the linear span of $Me_{A}M)$)

\item [c.~\cite{Po1,PiPo}] There exists  a semifinite trace  $\Tr$   on $\langle M,e_{A}\rangle$ determined by the equation $%
\Tr(xe_{A}y)=\tau (xy)$ for all $ x,y\in M$.

\item [d.~\cite{Po1,PiPo}] ($Pull$ $down$ $identity$) We have $e_{A}\langle M,e_{A}\rangle=\overline{e_{A}M}^{w}$ and $\langle M,e_{A}\rangle  e_{A}=\overline{Me_{A}}%
^{w}$. More precisely, there exists  an $M$--$M$ bimodule map $ \Phi :L^{1}(\langle M,e_{A}\rangle,\Tr)\rightarrow
L^{1}(M,\tau )$  which satisfies  the \emph{pull down identity}: if 
$V\in \langle M,e_{A}\rangle$, then  $ e_{A}V=e_{A}\Phi (e_{A}V)$.  Moreover,    $\Phi
(\langle M,e_{A}\rangle)  \subseteq L^{2}(M)$ .

\end{enumerate}

\end{proposition}
Next, we introduce the following definition:
\begin{definition}
A triple of von Neumann algebras $B \subseteq N \subseteq M$, where $M$ is a \II1 factor, has the relative weak asymptotic homomorphism property (relative WAHP) when \newline for all  $x_{1},x_{2},x_{3},...x_{n} \in M,$ and for every  $ \varepsilon
>0$ there exists  $u \in \unitary(B)$ such that: 
$$\left\Vert
E_{B}(x_{i}ux_{j})-E_{B}(E_{N}(x_{i})uE_{N}(x_{j}))\right\Vert
_{2}<\varepsilon $$ for $i,j=1,\dots,n$.
\end{definition}

For brevity,  we will often refer to the relative weak asymptotic homomorphism property
as relative WAHP. We were motivated to call this property {\em relative weak asymptotic homomorphism property} because it is a relative version for Robertson-Sinclair-Smith' s notion of {\em weak asymptotic homomorphism property} of a singular $MASA$. These facts will be more amply discussed in section 3.
In connection to the Definition above we record the following proposition:

\begin{proposition} Let M be  \II1  factor with $A\subseteq N\subseteq M$ two von Neumann subalgebras. The triple $A\subseteq N\subseteq M$ satisfy the relative weak asymptotic homomorphism property if and only if it satisfy the following property:

\medskip \noindent
  \quad For all $ x_{1},x_{2},\dots ,x_{n}\in M$, with $E_{N}(x_{i})=0$ and for every $\varepsilon >0$ there exists $u\in \mathcal{U}(A)$ such that  
$$||E_{A}(x_{i}ux_{j})||_{2}<\varepsilon $$
for all  $ i,j =  
1,2,...,n$
\medskip

\end{proposition}

\begin{proof}
The equivalence between relative weak asymptotic homomorphism property follows immediately if we use the identity:
$$E_{A}(xuy)-E_{A}(E_{N}(x)uE_{N}(y))=E_{A}((x-E_{N}(x))u(y-E_{N}(y)))$$\newline

\end{proof}

Next, we state the main theorem of this section:

\begin{theorem} \label{2.1}
Let M be \II1 factor and let $A\subset M$ be a MASA. 
 Then the triple $$A\subseteq  \mathcal{N}%
_{M}(A)^{\prime \prime} \subseteq M$$ satisfies 
the relative weak asymptotic homomorphism property.
\end{theorem}

The proof of this theorem will be a consequence of the next sequence of lemmas and propositions.
\begin{lemma}  Let M be \II1 factor and let $A\subset M$ be a MASA. Denote by $N:=\mathcal{N}%
_{M}(A)^{\prime \prime} $ and assume that the triple $A\subseteq  \mathcal{N}%
_{M}(A)^{\prime \prime}  \subseteq M$ does not satisfy relative WAHP.  
 Then, there exists 
a non-zero projection $f\in \mathcal{P}(A^{\prime
}\cap \langle M,e_{A}\rangle)$ such that $fe_{N}=0$, and $\Tr(f)<\infty $.
\end{lemma}

\begin{proof} 

The proof we present is essentially the proof of Corollary 2.3 in ~\cite{Po3} with very slight changes, and we reproduced it here for the sake of completeness.%

If the triple $A\subseteq  \mathcal{N}
_{M}(A)^{\prime \prime}  \subseteq M$ does not satisfy relative WAHP then, by Proposition 2.2, there exist 
$x_{1},x_{2},\dots ,x_{n}\in M$ such that: $E_N(x_{i})=0$ for $i=1,2,...,n$ 
and there exists $\varepsilon _{0}>0$ such that $$\Sigma _{i,j=1}^{n}||E_{A}(x_{i}ux_{j}^{\ast
})||_{2}^{2}\geq \varepsilon _{0}$$
  for all $ u\in \mathcal{U}(A).$
Using Proposition 2.1 c., a little computation shows that the last inequality is
equivalent to ($\alpha$)\quad $\Tr(bubu^{\ast })\geq \varepsilon _{0}$  for all $
u\in \mathcal{U}(A)$, where $b=\Sigma _{i=1}^{n}x_{i}^{\ast
}e_{A}x.$

Denote by $K(b)=\overline{co}^{w}\{ubu^{\ast
}|u\in \mathcal{U}(A)\}$ a  weak-operator-compact convex set in $\langle M,e_A\rangle$. If $a$ is the unique
$||$ $||_{2,\Tr}$-minimal element in $K(b)$ then:

$$0\leq a\leq 1,  \quad \Tr(a)\leq \Tr(b)\leq \infty,  \quad \text{ and  }
a\in A^{\prime
}\cap \langle M,e_{A}\rangle.$$

Also ($\alpha$) obviously implies that $\Tr(bx)\geq \varepsilon _{0}$  for all 
$x\in K(b)$. In particular $\Tr(ba)\geq
\varepsilon _{0}$ implies $a\neq 0$.
Since we are assuming that $E_{N}(x_{i})=0$ for all $i=1,2,...,n$, there are few more conditions, besides the details from the proof of Corollary 2.3 in ~\cite{Po3} that we need to check before we are be able to derive our conclusion.
Namely: 
$$be_{N}=\Sigma _{i=1}^{n}x_{i}^{\ast }e_{A}x_{i}e_{N}=\Sigma
_{i=1}^{n}x_{i}^{\ast }e_{A}e_{N}x_{i}e_{N}=\Sigma _{i=1}^{n}x_{i}^{\ast
}e_{A}E_{N}(x_{i})=0,$$

and because $e_{N}\in N^{\prime }\subset A^{\prime }$, we have
$$
\begin{aligned}
 &ubu^{\ast }e_{N}=0 \\
&\Rightarrow \co(ubu^{\ast })e_{N}=0\\
&\Rightarrow K(b)e_{N}=0  \quad \text{ and likewise }  e_{N}K(b)=0 \\
&\Rightarrow e_{N}a=ae_{N}=0 \text{ and}
 0\leq a\leq 1  \\
&\Rightarrow 0\leq a\leq 1-e_{N}.
\end{aligned}.
$$

Finally, by taking a suitable spectral projection of  $a$ in the algebra $%
(1-e_{N})A^{\prime }\cap \langle M,e_{A}\rangle(1-e_{N})$ we find a nonzero projection  $f\in
P(A' \cap \langle M,e_{A}\rangle)$ such that $fe_{N}=0$, and  $\Tr(f)<\infty $.

\end{proof}

\bigskip

\begin{lemma} Let $M$ be \II1 factor and let $A\subset M$ be a MASA. Suppose there exists $f\in A^{\prime }\cap \langle  M,e_{A} \rangle $, a nonzero projection that satisfies $\Tr(f)<\infty .
$. Then there exist nonzero projections $p_{i}\in A^{\prime }\cap \langle  M,e_{A} \rangle $ that are abelian 
in $\langle  M,e_{A} \rangle $ for all $i\in I$ and which further satisfy the following equation: $$f=\Sigma _{i\in I}p_{i}$$
\end{lemma}
\begin{proof} Denote by $\mathcal{A}:=A\vee \mathcal{J}{A}\mathcal{J}$
and consider the inclusions:%
\begin{equation*}
\mathcal{A}^{\prime }=A^{\prime }\cap \langle  M,e_{A} \rangle \subseteq \langle  M,e_{A} \rangle =(\mathcal{%
J}{A}\mathcal{J})^{\prime }
\end{equation*}
\begin{equation*}
fA^{\prime }\cap \langle  M,e_{A} \rangle f\subset f\langle  M,e_{A} \rangle f
\end{equation*}
\begin{equation*}
\mathcal{Z}{(\langle  M,e_{A} \rangle)}=\mathcal{%
J}{A}\mathcal{J}
\end{equation*}

First note that both algebras $A^{\prime }\cap \langle  M,e_{A} \rangle$ and $\langle  M,e_{A} \rangle$ are of type $I_{\infty }$. Since $f\in 
\mathcal{A}^{\prime }$ is a finite projection, both algebras $fA^{\prime }\cap \langle  M,e_{A} \rangle f$ and $f\langle  M,e_{A} \rangle f$
 are finite of type $I$ . Also note that the central support of $e_{A}$ in $\langle  M,e_{A} \rangle$ is equal to 1. ($z_{\langle  M,e_{A} \rangle }(e_{A})=1$)

By general theory, we have $f\langle  M,e_{A} \rangle f\cong \Sigma _{j\in J}$ $A_{j}%
\bar{\otimes}M_{n_{j}}(\mathbb{C})$. This implies that $\mathcal{A}f$
is an abelian subalgebra of $\Sigma _{j\in J}A_{j}\otimes M_{n_{j}}(\mathbb{C}).$

Consequently, $\mathcal{A}fz_{j}\subset A_{j}\bar{\otimes}M_{n_{j}}(\mathbb{C%
})$ (where $z_{j}\in \mathcal{J}{A}\mathcal{J}$) and there exists a MASA $B_{j}$,
such that $f\mathcal{A}z_{j}\subset B_{j}\subset A_{j}\bar{\otimes}M_{n_{j}}(%
\mathbb{C})$.

By Kadison's result ~\cite{K}, we have $uB_{j}u^{\ast }=A_{j}\otimes D_{n_{j}}(\mathbb{C%
})$, where $u\in \mathcal{U}(A_{j}\bar{\otimes}M_{n_{j}}\mathbb{C})$. That
 implies $fz_{j}=\Sigma _{i=1}^{n_{j}}p_{i}^{j}$ where $p_{i}^{j}\in (f%
\mathcal{A}z_{j})^{\prime }\cap fz_{j}Aj\otimes M_{n_{j}}(\mathbb{C})z_{j}f.$
But this shows that $p_{i}^{j}$ actually belongs to $\mathcal{A}^{\prime }fz_{j}$
and, moreover, is abelian in $\langle  M,e_{A} \rangle $  (i.e., $%
p_{i}^{j}\langle  M,e_{A} \rangle p_{i}^{j}=z_{j}p_{i}^{j}f\langle  M,e_{A} \rangle fp_{i}^{j}z_{j}=A_{j}%
\otimes \mathbb{C}p_{i}^{j}$ abelian algebra).

But this implies $f=\Sigma fz^{j}=\Sigma _{j,i\in \overline{1,n_{j}}}$ $%
p_{i}^{j}$ which completes the proof of this step. In particular, we have: if $f$ is a nonzero projection in $A^{\prime }\cap \langle  M,e_{A} \rangle $ with $\Tr(f)<\infty $, then there exists a nonzero
projection $p$ in $A^{\prime }\cap \langle  M,e_{A} \rangle $  which is abelian in $%
\langle  M,e_{A} \rangle $ and satisfies $p$ $\leq f$.
\end{proof}

\begin{proposition}\label{2.10} \mbox{} Let $A$ be a masa in the $II_1$ factor $M$ and let $N:=\mathcal{N}_M(A)^{\prime \prime}$. If $p$ is a non-zero projection in  $%
A^{\prime }\cap \langle  M,e_{A} \rangle $ with $p\preccurlyeq e_{A}$ as projections in $\langle  M,e_{A} \rangle$, then exists a non-zero projection $q$ with $q\leq p$ and $q\leq e_{N}$.
\end{proposition}

\begin{proof} Since  $p\preccurlyeq e_{A}$, let $W\in \langle  M,e_{A} \rangle $ be a partial isometry such
that $p=W^{\ast }W$, $WW^{\ast }\leq e_{A}$, $Wp=W$, $e_{A}W=W.$ \newline 

The property $e_A \langle M,e_A \rangle e_A=Ae_A$ (Proposition \ref{2.8} c) implies there exists an well defined function $\phi :A\rightarrow A$ given by the equation:  $$e_{A}WaW^{\ast }e_{A}=\phi (a)e_{A}.$$

 We remark that $\phi$ is a $*$-homomorphism. Indeed, it satisfies the following:

- $e_{A}W(a+b)W^{\ast }e_{A}=\phi (a+b)e_{A}$ and $e_{A}WaW^{\ast
}e_{A}+e_{A}WbW^{\ast }e_{A}=(\phi (a)+\phi (b))e_{A}\Rightarrow$ \quad $\phi
(a+b)=\phi (a)+\phi (b)$ by Proposition \ref{2.8} c.

- $\phi (a^{\ast })e_{A}=e_{A}Wa^{\ast }W^{\ast }e_{A}=(e_{A}WaW^{\ast
}e_{A})^{\ast }=\phi (a)^{\ast }e_{A}\Rightarrow \phi (a^{\ast })=\phi
(a)^{\ast }$.

- $\phi (ab)e_{A}=e_{A}WabW^{\ast }e_{A}=e_{A}WW^{\ast }e_{A}WabWe_{A}=e_{A}WaW^{\ast }e_{A}WbWe_{A}=$\quad $\phi (a)e_{A}\phi (b)e_{A}=\phi
(a)\phi (b)e_{A}\Rightarrow \phi (ab)=\phi (a)\phi (b)$ by Proposition \ref{2.8} c. again. Notice we also used here that $W^{\ast }e_{A}W\in A^{\prime }.$

Since $p=W^{\ast }e_{A}W\in A^{\prime }\cap
\langle  M,e_{A} \rangle $ we have $ap=pa$ far all $a\in A$. This implies $W^{\ast }e_{A}Wa=aW^{\ast }e_{A}W$ 
 and because 
$e_{A}W=W$ is an isometry we obtain $Wa=e_{A}WaW^{\ast}e_{A}W$, which we rewrite as $$Wa=\phi(a)W$$   
 for all $a\in A$.
 
 For a more detailed account on $*$-homomorphism $\phi$ we send the reader to ~\cite{Po3,PoSiSm}.

Since span$Me_AM$ is weakly dense $*$-algebra in $\langle  M,e_{A} \rangle$, using Kaplansky density theorem, there exists $(z_n)_n\in$ span$Me_AM $ such that $z_n\rightarrow W^{\ast }$ in $so$-topology. But $W^{\ast }=W^{\ast }e_A$ implies that  $z_ne_A\rightarrow W^{\ast }e_A=W^{\ast }$. Using $e_AMe_A=Ae_A$ we see that $z_ne_A=y_ne_A$ with $y_n\in M$. Next, denote by $\eta=\mathcal{J}W^{\ast }\hat{1}\in L^2(M,\tau)$ where $\hat{1}$ is the canonical cyclic trace vector for the left regular representation associated with $\tau$.

Also by $l_{\eta}$ we mean the left multiplication operator by $\eta$ defined on  $M\hat{1}$. It is well known this operator is closable and we denote by $L_\eta:=\overline{l_{\eta}}$ its closure. We record that $L_\eta$ is a closed densely defined operator affiliated with $M$. \newline
Next equation establishes the relation between $W$ and $L_\eta$, which is known in the literature as the pull-down identity (see ~\cite{Po1,PiPo}:\newline
$$e_{A}Wx\hat{1}=e_AL_\eta(x\hat{1})\qquad (2.6.1.)$$ for all $x\in M$.
 
 To verify this we observe that it is enough to check the following:\newline
 $\langle e_AL_\eta(x\hat{1}),y\hat{1}\rangle=\langle L_\eta(x\hat{1}),E_A(y)\hat{1}\rangle$ \newline $=\langle x\hat{1},L_\mathcal{J}\eta(E_A(y)\hat{1})\rangle=\langle x\hat{1},\mathcal{J}E_A(y^{\ast})\mathcal{J}W^{\ast }\hat{1})\rangle$ \newline $=\lim_{n} \langle x\hat{1},\mathcal{J}E_A(y^{\ast})\mathcal{J}z_n\hat{1})\rangle=\lim_{n} \langle x\hat{1},z_nE_A(y)\hat{1})\rangle$ \newline $=\lim_{n} \langle x\hat{1},z_ne_AE_A(y)\hat{1})\rangle=\lim_{n}\langle x\hat{1},z_ne_AE_A(y)\hat{1})\rangle$ \newline $=\langle x\hat{1},W^{\ast}E_A(y)\hat{1})\rangle=\langle Wx\hat{1},E_A(y)\hat{1})\rangle=\langle e_AWx\hat{1},y\hat{1})\rangle$  for every $x,y\in M$.\newline

 From now on, whenever two unbounded operators $S$ and $T$ defined on $L^2(M,\tau)$ agree on $M\hat{1}$, we write  $S\dot{=}T$. Consequently, we can rewrite the equation (2.6.1.) as:
  
  $$e_AW\dot{=}e_AL_{\eta}.$$
  
  Further, we can also check the following:\newline
  $L_{\eta}a(x\hat{1})=L_{\eta}((ax)\hat{1})$\newline$=\mathcal{J}(ax)^*\mathcal{J}\eta=\mathcal{J}x^*a^*\mathcal{J}\mathcal{J}W^*\hat{1}$\newline$=\mathcal{J}x^*a^*W^*\hat{1}=\mathcal{J}x^*W^*\phi(a^*)\hat{1}$\newline$=\mathcal{J}x^*W^*\mathcal{J}\phi(a)\hat{1}=\phi(a)\mathcal{J}x^*\mathcal{J}\mathcal{J}W^*\mathcal{J}\hat{1}$\newline$=\phi(a)\mathcal{J}x^*\mathcal{J}{\eta}=\phi(a)L_{\eta}(x\hat{1})$ for every $a\in A$ , $x\in M$.
  In other words we have proved that:
 
  $$L_{\eta}a\dot{=}\phi(a)L_{\eta}\qquad(2.6.2.)$$ for all $a\in A$.

  Let  $L_{\eta}=wT$ be the polar decomposition of $L_{\eta}$ where $w\in M$ is the partial isometry mapping the closure of the range of $T$ to the closure of the range of $L_{\eta}$ and $T=|L_{\eta}|$ is the absolute value of $L_{\eta}$.
  
  The relation (2.6.2.) becomes $$wTa\dot{=}\phi(a)wT\qquad(2.6.3.)^{\prime}$$ which by the same argument like in the Lemma 5.1 from ~\cite{PoSiSm} implies that
  
   $$w^*wa=w^*\phi(a)w\qquad (2.6.3.)$$ for all $a\in A$.\newline
   
   In particular this equation carries the fact that $w^*w\in A$ and, moreover, it can be proved that $ww^*\in \phi(A)^{\prime}\cap M$.
   
   To see this let $f\in \mathcal{P}(A)$ be an arbitrary fixed projection. Then, $w^*\phi(f)ww^*\phi(f)w=w^*wfw^*wf=w^*wf^2=w^*wf=w^*\phi(f)w$, which further implies that $ww^*\phi(f)(1-ww^*)=0$ and consequently $ww^*\phi(f)=\phi(f)ww^*$. Since this last equation holds true for any projection $f\in A$ our conclusion follows.
   So the equation (2.6.3.) is actually equivalent to $$wa=\phi(a)w\qquad(2.6.4.)$$ for all $a\in A$.
   
   Also, if we further combine (2.6.4.) with the fact that $w^*wT\dot{=}T$ then the equation $(2.6.3.)^{\prime}$ implies: 
   
   $$aT\dot{=}Ta\qquad (2.6.5)$$ for all $a\in A$.

Now consider $\mathcal{D}_1:=\mathcal{D}(\overline{aT})\cap \mathcal{D}(\overline{Ta})\cap \mathcal{D}(T)$, where by  $\mathcal{D}(\overline{aT})$ we denoted the domain of the (closed) operator $\overline{aT}$ and so on....
         
 Because $M\hat{1}\subset \mathcal{D}(T)\subset\mathcal{D}(aT)\subset \mathcal{D}(\overline{aT})$
and  $M\hat{1}\subset \mathcal{D}(Ta)\subset\mathcal{D}(\overline{Ta})$, we have that $M\hat{1}\subset \mathcal{D}_1$.
 Now by Lemma 16.4.3 in ~\cite{MvN} we have that $\mathcal{D}_1$ is essentially dense and everywhere dense in $L^2(M,\tau)$
 and moreover we prove below that the equality (2.6.5.) actually holds true on $\mathcal{D}_1$.
 
  If $\zeta\in \mathcal{D}_1$ then there exists 
$m_n\hat{1} \in M\hat{1}$ such that  $m_n\hat{1}\rightarrow\zeta$. Next, the equation (2.6.5.) implies that $\langle m_n\hat{1},\overline{\overline{aT}-\overline{Ta}}(m_n\hat{1})\rangle=\langle m_n\hat{1},aT-Ta(m_n\hat{1})\rangle=\langle m_n\hat{1},0\rangle \rightarrow \langle \zeta,0 \rangle$ which futher gives that $\overline{\overline{aT}-\overline{Ta}}(\zeta)=0$ and so $\overline{aT}(\zeta)=\overline{Ta}(\zeta)$ for every $\zeta\in \mathcal{D}_1$. Since $Ta=\overline{Ta}$ and $\mathcal{D}(T)\subset\mathcal{D}(aT)$ we conclude that 
 
$$aT(\zeta)\dot{=}Ta(\zeta)\qquad(2.6.6.)$$ for all $\zeta\in \mathcal{D}_1$ and all $a\in A$.

Now consider $\mathcal{D}_2:=\mathcal{D}_1\cap \mathcal{D}(\overline{T^2})\subset \mathcal{D}_1$. By Lemma 16.4.3 ~\cite{MvN} again we have that $\mathcal{D}_2$ is essentially dense and everywhere dense in $L^2(M,\tau)$. 

Next we verify that ${T^2+a^2}_{\mid_{\mathcal{D}_2}}$ is essentially selfadjoint. On one hand we have
$\langle T^2+a^2\zeta,\psi \rangle=\langle\zeta,\overline{T^2}+a^2\psi \rangle$ and hence 

${(T^2 +a^2}_{\mid_{\mathcal{D}_2}})^{*}\subset \overline{{T^2}+a^2}_\mid{_{\mathcal{D}(\overline{T^2})}}$ which is obviously a closed operator. Using the uniqness of the extension from Lemma 16.4.2 ~\cite{MvN} we see that

$({T^2+a^2}_{\mid_{\mathcal{D}_2}})^*=\overline{{T^2}+a^2}_{\mid_{\mathcal{D}(\overline{T^2})}}$. On the other hand, 
$\overline{{T^2+a^2}_{\mid_{\mathcal{D}_2}}}\subset \overline{{\overline{T^2}+a^2}_{\mid_{\mathcal{D}(\overline{T^2})}}}={\overline{T^2}+a^2}_{\mid_{\mathcal{D}(\overline{T^2})}}$ and by uniqueness again we have $\overline{{T^2+a^2}_{\mid_{\mathcal{D}_2}}}={\overline{T^2}+a^2}_{\mid_{\mathcal{D}(\overline{T^2})}}$.\newline

In conclusion ${T^2+a^2}_{\mid_{\mathcal{D}_2}}$ is essentially selfadjoint. Combining this with the fact that $T$ and $a$ commutes on $\mathcal{D}_2$ (see $(2.6.6.)$), by Corollary 9.2 in ~\cite{Ne} we obtain that $T$ and $a$ strongly commute which means their spectral scales commutes. So the spectral scale of $T$ belongs to $A^{\prime}\cap M=A$.

 To this end we prove the following:

 \emph{Claim:  There exists $m\in M$ such that $m^*e_Am$ is a nonzero projection in $A^{\prime}\cap\langle M,e_A\rangle$ that satisfies $m^*e_Am\leq p$.} \newline

To show this let $f$ be a spectral projection of $T$ such that $0\neq fT=Tf\in M$. By (2.6.1.) we have  $Wf(x\hat{1})=e_AW(fx\hat{1})=e_AL_{\eta}(fx\hat{1})=e_AwT(fx\hat{1})=e_Aw(Tf)(x\hat{1})$. Since $Tf$ is a bounded operator that belongs to $M$ we have obtained that
$Wf=e_Am$ for some $m\in M$.

Also by the choice of $f$ we have that $0\neq Wf$ so we can verify the following:
\indent$0\neq m^*e_Am=fW^*Wf=W^*Wf^2=W^*WfW^*W\leq (W^*W)^2=p$\newline \indent$m^*e_Am=fW^*Wf=(W^*Wf)^2=(me_Am^*)^2$.\newline
\indent$m^*e_Am=fW^*Wf=W^*Wf\in A^{\prime}\cap\langle M,e_A\rangle$. \newline
This finishes the proof of the \emph{Claim}.\newline

Since $m^*e_Am\in A^{\prime}$ we have 
$m^*e_Ama=am^*e_Am$ for all $a\in A$, which is equivalent to 
$$ma=E_{A}(mam^*)m\qquad (2.6.7.)$$ for all $a\in A$.

Because relation (2.6.7.) holds for every $a\in A$ by considering the "stared" version of it we obtain $m^*m\in A^{\prime}\cap M=A$ and hence $|m|\in A$.

 Next, denote by $h=\chi_{(\infty,1)}(|m|)\in A$ the spectral projection of the element $|m|$ corresponding to the interval $(\infty,1)$.
 To this end we split the proof of this proposition in two cases:\newline

 \emph{CASE I $h\neq 0$.} 
 
 Relation (2.6.7.) implies  that $mha=E_{A}(mhahm^*)mh\quad (2.6.8)$ for all $a\in A$.
 Consider $m_1=mh$ and using the spectral properties of $h$ we see $ m_1^*m_1=hm^*mh\leq h\leq 1$ which is  equivalent to
 $m_1m_1^*\leq 1$.
 Also, by plugging in $a=1$ in equation (2.6.8.)  we have $m_1= E_{A}(m_1m_1^*)m_1$ which in particular implies that $E_{A}(m_1m_1^*)\in\mathcal{P}(A)$.

Next, notice  
$m_1m_1^*=E_{A}(m_1m_1^*)m_1m_1^*E_{A}(m_1m_1^*)\leq (E_{A}(m_1m_1^*))^2=E_{A}(m_1m_1^*)$ and by the $\tau$-invariance of $E_{A}$ and the faithfullness of $\tau$ we conclude $m_1m_1^*=E_{A}(m_1m_1^*)\quad(2.6.9.)$.

Also we remark that the equation (2.6.8.) implies that $m_1^*m_1\in A$ and this together with equations (2.6.9.) help us to conclude that $m_1\in \mathcal{GN}_M(A)$. Using the structure of $\mathcal{GN}_M(A)$(~\cite{D}) there exists $u\in \mathcal{N}_M(A)$, $e\in A$ such that $m_1=ue$ and we can verify that $m_1^*e_Am_1=P_{m_1^*A}=P_{eu^*A}\leq P_N=e_N$.

But $h\neq 0$ implies $m_1^*e_Am_1\neq 0 $ and also we can check that\newline $m_1^*e_Am_1=hm^*e_Amh=m^*e_Amhm^*e_Am\leq m^*e_Am\leq p$.
Hence in this case the proof of the proposition is finished.\newline

\emph{CASE II $h=0$.}

 If $h=0$ , then $\sigma (|m|)\subset \lbrack 1,\infty )$ so in particular $|m|$ is invertible and \newline $1\leq m^*m\quad(2.6.10.)$. By taking the polar decomposition of $m=v|m|$ we have that \newline $v^*v=supp(|m|)=1$ which implies  $v\in \mathcal{U}(M)$ because $M$ is a finite factor. 
 
 Moreover, the  relation $(2.6.7.)$ becames
$v|m|a=E_{A}(mam^*)v|m|$ and so $va|m|=E_{A}(mam)^*v|m|$ for all $a\in A$.
By multiplying on the right by $|m|^{-1}$ we get $va=E_{A}(mam^*)v$ which we rewrite as $vav^*=E_{A}(mam^*)$ for all $a\in A$. This last equation implies $vAv^*\subseteq A$ which together with $v\in \mathcal{U}(M)$ and $A$ is a masa in $M$ further implies that 
$v\in \mathcal{N}_M(A)$. 

From this we notice that $0\neq v^*e_Av=P_{v^*A}\leq e_N\quad(2.6.11.)$ and $v^*e_Av\in A^{\prime}\cap\langle M,e_A \rangle$.

On the other hand using $(2.6.10.)$ we can check the following 

$0\neq v^*e_Av\leq v^*e_Avm^*mv^*e_Av=v^*e_Av|m|^{2}v^*e_Av=|m|v^*e_Av|m|=m^*e_Am^*\leq p\quad (2.6.12.)$

Relations $(2.6.11)$ and $(2.6.12.)$ finish the proof in this case.

\end{proof}

 We end this section by presenting \emph{the proof of the Theorem 2.3}:

\begin{proof}: We will proceed by contradiction. Let suppose that the triple does not satisfy the relative WAHP.
By Lemma 2.4 there exists a nonzero projection $f\in \mathcal{P}(A^{\prime }\cap \langle  M,e_{A} \rangle )\ $with $%
\Tr(f)<\infty$ and $fe_N=0\quad (\ddagger) $. Moreover, using Lemma 2.5 there exists a nonzero projection $p\leq f$ which lies in  $A^{\prime }\cap \langle  M,e_{A} \rangle$ and is abelian $\langle  M,e_{A} \rangle$. But we obviously have $z_{\langle  M,e_{A} \rangle }(p)\leq 1=z_{\langle  M,e_{A} \rangle }(e_{A})$ which further implies that $p\preccurlyeq e_{A}$. By Proposition 2.6 there exists a non-zero projection $q$ such that $q\leq p\leq f$ and $q\leq
e_{N}$, which is in contradiction with $(\ddagger) $. In conclusion the triple $A\subseteq \mathcal{N}_M(A)^{\prime \prime}\subseteq M$ must satisfy the relative WAHP. 
\end{proof}

\section{ Applications}

In this section we present several immediate applications of Theorem \ref{2.1}. \newline

The first result of the section underlines the fact that the \emph{relative weak asymptotic homomorphism property} (relative WAHP) for a triple of algebras is a suitable tool to control the normalizing algebra of a given subalgebra.
This idea was exploited before in \cite{RoSiSm} (see Lemma 2.1), where the authors proved that the \emph{weak asymptotic homomorphism property} (WAHP) of a MASA implies the strong singularity of that MASA. The same argument can be used to prove the following:

\begin{proposition} \label{3.1}

Let M be a $II_{1}$ factor and $A\subseteq N\subseteq M$ two von Neumann subalgebras. If the triple $A\subseteq N\subseteq M$ satisfies the relative weak asymptotic homomorphism property then the following inequality holds for every $u\in \mathcal{U}(M)$:
 
\begin{equation*}
||E_{uAu^{\ast }}-E_{A}||_{\infty ,2}\geq ||u-E_{N}(u)||_{2} 
\end{equation*}

\end{proposition}

\begin{proof} 
If we apply relative WAHP for the set $\{u,u^{\ast }\}$ and $\varepsilon >0$
arbitrary but fixed, then there is $a_{\varepsilon }\in \mathcal{U}(A)$ such
that: 
\begin{equation*} 
(\gamma)\qquad||E_{A}(u^{\ast }a_{\varepsilon }u)-E_{A}(E_{N}(u^{\ast
})a_{\varepsilon }E_{N}(u))||_{2}<\varepsilon. 
\end{equation*}%

But,\newline 
$||E_{uAu^{\ast }}-E_{A}||_{\infty ,2}^{2}\geq ||E_{uAu^{\ast
}}(a_{\varepsilon })-E_{A}(a_{\varepsilon })||_{2}^{2}$

\qquad \qquad \qquad \qquad$=1-||E_{A}(u^{\ast }a_{\varepsilon }u)||_{2}^{2}$ (because $a_{\varepsilon }\in \mathcal{U}(A)$)

\qquad \qquad \qquad \qquad$\geq 1-(||E_{A}(E_{N}(u^{\ast })a_{\varepsilon }E_{N}(u)||_{2}+\varepsilon
)^{2}$ (by$(\gamma)$)

\qquad \qquad \qquad \qquad$\geq 1-(||E_{N}(u^{\ast }a_{\varepsilon }E_{N}(u)||_{2}+\varepsilon )^{2}$ (since $e_{A}$ is a projection in $B(L^{2}(M)))$

\qquad \qquad \qquad \qquad$\geq 1-(||E_{N}(u)||_{2}+\varepsilon )^{2}$ ( $E_{A}$ is a Schwartz map, $a_{\varepsilon }\in \mathcal{U}(A)$, $u\in \mathcal{U}(M))$

\qquad \qquad \qquad \qquad$\geq ||u -E_{N}(u)||_{2}^{2}-\varepsilon (\varepsilon +2)$\newline

Note this is true for any $\varepsilon >0$ so by taking $\varepsilon
\rightarrow 0$ we obtain the desired result. 
\end{proof}
\begin{corollary} \label{3.2} \mbox{}

Let $M$ be a $II_{1}$ factor with $A$ and $N$ two von Neumann subalgebras. If we assume that the triple $A\subseteq N\subseteq M$ satisfies the relative WAHP,
 then we have ${\mathcal{N}%
_{M}(A)}^{\prime \prime}\subseteq N$.

\end{corollary}
\begin{proof}
The proof is an obvious consequence of the previous proposition.
\end{proof}

\begin{theorem} \label{3.2} \mbox{}

 (A generalization of Sinclair-Smith inequality) Let $M$ be a  $II_{1}$ factor with $A\subset M$ a MASA. Then,

$$||E_{uAu^{\ast }}-E_{A}||_{\infty ,2}\geq ||u-E_{\mathcal{N}_{M}(A)}(u)||_{2}, \:{ for}\: {all }\: u\in \mathcal{U}(M).$$

\end{theorem}
  
\begin{proof}By Theorem \ref{2.1} the triple $A\subseteq \mathcal{N}_M(A)^{\prime\prime}\subseteq M $ satisfies the relative WAHP so the statement follows from Proposition \ref{3.1}.
\end{proof}

At this point we would like to mention that even though the previous theorem is a generalization of strong singularity concept for a MASA in a \II1 factor, a more general version (with absolute constant 1) for Theorem 6.2 in ~\cite{PoSiSm} remains open.

\medskip
In the last part of this section we will present a series of estimates of the normalizing algebra of a MASA in the situations of tensor products and cross products by discrete groups. These estimates  heavily rely on relative  WAHP for certain triples of algebras. Before starting we state an alternative description of relative WAHP for a triple of albgebras which is more convenient to use in our future computations.

\begin{remark}

Let M be a $II_{1}$ factor and $A\subseteq N\subseteq M$ two von Neumann subalgebras. Let $\mathcal{X}\subset M$ such that span$\mathcal{X}$ is a *-subalgebra which is weakly dense in $M$. The triple $A\subseteq N\subseteq M$ satisfies the \emph{relative weak asymptotic homomorphism property} if and only if:

For all   $x_{1},x_{2},...x_{n} \in \mathcal{X} ,$ and for every  $ \varepsilon
>0$ there exists  $u \in \unitary(A)$ such that 
$$\left\Vert
E_{A}(x_{i}ux_{j})-E_{A}(E_{N}(x_{i})uE_{N}(x_{j}))\right\Vert
_{2}<\varepsilon $$ for $i,j=1,\dots,n.$

\end{remark}

\begin{proposition} \label{3.3}
\ignore{
Let $M_2\vee M_1,M_1,M_2$ $\ II_{1}$ factors with $[M_1,M_2]=0$, $M_2\supset A_{2}$ MASA, $%
M_1\supset A_{1}$ MASA, $N_{1}:=\mathcal{N}_{M_1}(A_{1})^{\prime \prime }$, $%
N_{2}:=\mathcal{N}_{M_2}(A_{2})^{\prime \prime }$. Consider the following
diagram   $%
\begin{array}[t]{ccccccc}
M\vee N & \supset  & M_2 & \supset  & N_{2} & \supset  & A_{2} \\ 
\cup  &  &  &  &  &  &  \\ 
M_1 &  &  &  &  &  &  \\ 
\cup  &  &  &  &  &  & \cup  \\ 
N_{1} &  &  &  &  &  &  \\ 
\cup  &  &  &  &  &  &  \\ 
A_{1} &  &  & \supset  &  &  & B%
\end{array}%
$
}
Let   $M_1$ and $M_2$ be $\ II_{1}$ factors.    For  $i = 1, 2$,  let $A_i$ be a MASA in $M_i$.
Then, $$\mathcal{N}_{M_1\overline{\otimes} M_2}(A_{1}\overline{\otimes} A_{2})^{\prime\prime}=\mathcal{N}%
_{M_2}(A_{2})^{\prime\prime}\overline{\otimes} \mathcal{N}_{M_1}(A_{1})^{\prime\prime}.$$%
In particular if $A_i$ is a singular MASA in $M_i$ for  $i = 1, 2$ then $A_{1}\overline{\otimes} A_{2}\subset M_1\overline{\otimes} M_2$ is a singular MASA. 
\end{proposition}

\begin{proof}  For each $i = 1, 2$, we denote by  $N_{i}:=\mathcal{N}_{M_i}(A_{i})^{\prime\prime}$. We only need to prove $$\mathcal{N}_{M_1\overline{\otimes} M_2}(A_{1}\overline{\otimes}
A_{2})^{\prime\prime}\subseteq \mathcal{N}_{M_2}(A_{2})^{\prime\prime}\overline{\otimes}
\mathcal{N}_{M_1}(A_{1})^{\prime\prime},$$ the other containement being trivial.

By Corollary \ref{3.2}\, to show this would be enough to prove that the triple: 

$$A_1\overline{\otimes} A_2\subseteq N_1\overline{\otimes} N_2\subseteq M_1\overline{\otimes} M_2$$ satisfies the relative WAHP.

Further, since $\mathcal{X}=$span \{$ x\otimes y \mid\quad x\in M_1, y\in  M_2 $\} is $\parallel\quad\parallel_2$-dense in $M_1\overline{\otimes} M_2$ the triple $A_1\overline{\otimes} A_2\subseteq N_1\overline{\otimes} N_2\subseteq M_1\overline{\otimes} M_2$ satisfies the relative WAHP iff we have the following: \newline
$(*)$\quad For all $x_{1}\otimes y_1,\dots ,x_{n}\otimes y_n \in
\mathcal{X}$ and every $\varepsilon >0$ exists $a\in\mathcal{U}(A_{1}\overline{\otimes} A_{2})$ such
that: $$||E_{A_{1}\overline{\otimes} A_{2}}((x_{i}\otimes y_i)a(x_{j}\otimes y_j)-E_{A_{1}\overline{\otimes}
A_{2}}(x_{i}\otimes y_iaE_{N_{1}\overline{\otimes} N_{2}}(x_{j}\otimes y_j)||_{2}<\varepsilon $$ \ for all $i,j=1,\dots ,n.$\newline

 By Theorem \ref{2.1} we have that the triples $A_{1}\subseteq N_{1}\subseteq M_{1}$ and $A_2\subseteq N_{2}\subseteq M_{2}$ satisfies the relative WAHP, so there exists $a_{1}\in\mathcal{U} (A_{1})$, $a_{2}\in\mathcal{U}(A_{2})$ which satisfies the following inequalities:

$$||E_{A_{1}}(x_{i}a_{1}(x_{j}-E_{N_{1}}(x_{j})))||_{2}<\frac{\varepsilon }{%
2\underset{i=1,\dots ,n}{max}||y_{i}||^{2}}$$ \ for all $i,j=1,\dots ,n$
and 

$$||E_{A_{}}(y_{i}a_{2}(y_{j}-E_{N_{2}}(y_{j})))||_{2}<\frac{\varepsilon }{%
2\underset{i=1,\dots,n}{max}||y_{i}||^{2}}$$
for all $i,j=1,\dots ,n$.
 
 Finally, we evaluate:
$$
\begin{aligned}
&||E_{A_{1}\overline{\otimes} A_{2}}((x_{i}\otimes y_{i})(a_{1}\otimes a_{2})(x_{j}\otimes y_{j}-E_{N_{1}\overline{\otimes}
N_{2}}(x_{j}\otimes y_{j})))||_{2} \\
&=||E_{A_{1}}(x_{i}a_{1}(x_{j}-E_{N_{1}}(x_{j}))\otimes E_{A_{2}}(y_{i}a_{2}y_{j})+E_{A_{1}}(x_{i}a_{1}E_{N_{1}}(x_{j}))\otimes E_{A_{2}}(y_{i}a_{2}(y_{j}-E_{N_{2}}(y_{j}))||_{2} \\
&\leq \underset{i=1,\dots,n}{max}
||y_{i}||^{2}||E_{A_{1}}(x_{i}a_{1}(x_{j}-E_{N_{1}}(x_{j})||_{2}+\underset{i=1,\dots,n
}{max}||x_{i}||^{2}||E_{A_{2}}(y_{i}a_{2}(y_{j}-E_{N_{2}}(y_{j}))||_{2} \\
&<\varepsilon
\end{aligned}
$$
This completes the proof of $(*)$ and Proposition 3.5.

\end{proof}

\begin{corollary}\label{3.4}
If $M_i$ is a   \II1 factor  and  $A_i \subseteq M_i$ is a MASA for every $i=1,\dots ,k$ 
 then  
we have: 
$$\overline{\otimes }_{i=1}^{k}   \mathcal{N}%
_{M_{i}}(A_{i})^{\prime\prime}=\mathcal{N}_{\overline{\otimes }%
_{i=1}^{k}M_{i}}(\overline{\otimes}A_{i})^{\prime \prime}  \quad  (\ast\ast )
$$
\end{corollary}
\begin{proof}
The proof follows from Proposition \ref{3.3} by induction.
\end{proof}
\medskip
At this point it is natural to investigate if $(\ast\ast )$ holds true for
infinite tensor products. As expected, the answer is yes but to be able to prove this we first need to 
 analyze the behavior of the normalizing algebra with respect to the inductive limit.

\bigskip

\begin{proposition}\label{3.5}
\ignore{
Let $M$ $II_{1}$ factor and suppose the following holds true:\newline
$%
\begin{array}[t]{ccccccc}
M & \supset  & \overline{\cup _{n}N_{n}} & \supset  & \bigvee_{n}P_{n} & 
\supset  & \overline{\cup A_{n}}^{W} \\ 
&  & \cup 
 &  &  &  & \cup  \\ 
&  & \vdots  &  &  &  & \vdots  \\ 
&  & \cup  &  &  &  & \cup  \\ 
&  & N_{n+1} & \supset  & P_{n+1} & \supset  & A_{n+1} \\ 
&  & \cup  &  &  &  & \cup  \\ 
&  & N_{n} & \supset  & P_{n} & \supset  & A_{n} \\ 
&  &  &  &  &  & 
\end{array}%
$ with $A_{n}\subset N_{n}$ MASA, $N_{n}$ is a $II_{1}$factor, $%
P=\bigvee_{n}P_{n}$, $N=\overline{\cup _{n}N_{n}}^{W}$, $A=\overline{\cup
_{n}A_{n}}^{W}$,$P_{n}=\mathcal{N}_{N_{n}}(A_{n})^{\prime \prime }$ 
}

Let   $M_n$,   $n \in \mathbb N$,    be an increasing sequence of  $II_1$ factors that are contained in a larger $II_{1}$ factor $Q$.
Let $A_n \subseteq M_n$ be a MASA,  and suppose that  $A_n \subseteq A_{n+1}$  and moreover     
$$ 
\begin{array}[t]{ccc}
M_{n+1} & \supset  & A_{n+1} \\ 
\cup  &  & \cup  \\ 
M_{n} & \supset  & A_{n} \\ 
&  & 
\end{array}
$$
 is a commuting square  for all $n \in \mathbb N$. Denote by $P_n = \mathcal{N}_{M_{n}}(A_{n})^{\prime\prime}$, by $M=\overline{\cup _{n}M_{n}}^{w}$ and by $A=\overline{\cup
_{n}A_{n}}^{w}$. Then $A\subset N$ is MASA and $\mathcal{N}_{M}(A)^{\prime\prime}\subseteq
\bigvee_{n}P_{n}$.
\end{proposition}

\begin{proof}:   It is helpful to keep the following diagram in mind:
$$%
\begin{array}[t]{ccccccc}
 &  & M & \supset  & \bigvee_{n}P_{n} & 
\supset  & A \\ 
&  & \cup 
 &  &  &  & \cup  \\ 
&  & \vdots  &  &  &  & \vdots  \\ 
&  & \cup  &  &  &  & \cup  \\ 
&  & M_{n+1} & \supset  & \overset{n+1}{\underset{s=1}{\bigvee}}{P_s} & \supset  & A_{n+1} \\ 
&  & \cup  &  &  &  & \cup  \\ 
&  & M_{n} & \supset  & \overset{n}{\underset{s=1}{\bigvee}}{P_s} & \supset  & A_{n} \\ 
&  &  &  &  &  & 
\end{array}%
$$

We will only prove the second statement, the first one being nothing but Proposition 5.2.2 in ~\cite{Po5}.     Write $P = \bigvee_{n}P_{n}$.  In order to prove that 
 $\mathcal{N}_{M}(A)^{\prime \prime}\subseteq P$,
 it is enough to show that the triple $A \subseteq P \subseteq M$   satisfies the relative WAHP, which by Remark 3.4  reduces to showing the following:  \newline  
 For every $%
\varepsilon >0$,  for each $n \in N$ and for all $x_{1},\dots x_{k}\in M_n$ there exists $a\in \mathcal{U}(A)$ such that $$%
||E_{A}(x_{i}ax_{j})-E_{A}(E_P(x_{i})aE_P(x_{j}))||_{2}<\varepsilon $$  for all $ i,j=1,...,k$.

Let $\varepsilon >0$, $n \in N$  and $x_{1},\dots, x_{k}\in M_n$ fixed. Obviously we have that there exists $l\in \mathbb{N}$ such that $n\leqslant l$ and $$||E_{\overset{l}{\underset{s=1}{\bigvee}}{P_s}}(x_{i})-E_{P}(x_{i}))||_{2}<\frac{\varepsilon }{%
2\underset{i=1,\dots,k}{max}||x_{i}||^{2}}$$ for  every $1 \le i \le k$.

By Theorem 2.3 we know that the triple $A_{n}\subseteq P_{n}\subseteq M_{n}$ verifies the relative WAHP,  which implies that the triple $A_{n}\subseteq P_{n}\subseteq M_{n}$ also satisfies the relative WAHP,  so we have that there exists  $a_{l}\in 
\mathcal{U}(A_{l})\subset \mathcal{U}(A)$ such that: $$(\beta)\qquad||E_{A_{l}}(x_{i}a_{l}x_{j})-E_{A_{l}}(E_{\overset{l}{\underset{s=1}{\bigvee}}{P_s}}(x_{i})a_{l}x_{j})||_{2}<%
\frac{\varepsilon}{2} $$ for all $i,j=1,...,k$. 

Proceeding like in the proof of Proposition 4.2.2. in ~\cite{GHJ} we have that the commuting square condition is preserved under inductive limit.

Consequently, we have:
$$
\begin{array}[]{ccc}
M & \supset  & A \\ 
\cup  &  & \cup  \\ 
M_{l} & \supset  & A_{l} \\ 
&  & 
\end{array}%
$$
is a  commuting square for all $l\in \mathbb{N}$.

The commuting square condition enable us to observe that the inequalities $(\beta)$
are equivalent to 
$$||E_{A}(x_{i}a_{l}x_{j})-E_{A}(E_{\overset{l}{\underset{s=1}{\bigvee}}{P_s}}(x_{i})a_{l}x_{j})||_{2}< \frac{\varepsilon}{2} $$ for all $i,j=1,...,k$.

Finally, we evaluate:

 $||E_{A}(x_{i}a_lx_{j})-E_{A}(E_P(x_{i})a_lE_P(x_{j}))||_{2}=||E_{A}(x_{i}a_lx_{j})-E_{A}(E_P(x_{i})a_lx_{j})||_{2} \leq$ \newline
 
 $||E_{A}(x_{i}a_lx_{j})-E_{A}(E_{\overset{l}{\underset{s=1}{\bigvee}}{P_s}}(x_{i})a_lx_{j})||_{2}+||E_{A}(E_{\overset{l}{\underset{s=1}{\bigvee}}{P_s}}(x_{i})a_lx_{j})-E_{A}(E_P(x_{i})a_lx_{j})||_{2}<$\newline $$ \frac{\varepsilon}{2} +  \underset{i=1,\dots,k}{max}||x_{i}||^{2} \frac{\varepsilon }{%
2\underset{i=1,\dots,k}{max}||x_{i}||^{2}}=\varepsilon$$ \newline and we are done.

\end{proof}
\medskip

\begin{corollary}\label{3.6}

Let $I$ be a countable set, $\{M_i\}_{i\in I}$ a collection of $II_1$ factors and  $\{A_i\}_{i\in I}$ a collection of abelian diffuse von Neumann algebras such that for every $i\in I$ we have that $A_i\subset M_i$ is a MASA. Then, $$\overline{\bigotimes}_{i\in I}(\mathcal{N}_{M_{i}}(A_{i})^{\prime\prime})= (\mathcal{N}_{\overline{\bigotimes}_{i\in I}M_{i}}(\overline{\otimes}%
_{i\in I} A_{i}))^{\prime\prime}.$$
\end{corollary}

\begin{proof} It follows immediately by applying Proposition \ref{3.5} and Corollary \ref{3.4} together with
the fact that \newline
$$%
\begin{array}[t]{ccc}
\overline{\bigotimes }_{i\in S_{n+1}}M_{i} & \supset  & \overline{\bigotimes 
}_{i\in S_{n+1}}A_{i} \\ 
\cup  &  & \cup  \\ 
\overline{\bigotimes }_{i\in S_{n}}M_{i} & \supset  & \overline{\bigotimes }%
_{i\in S_{n}}A_{i} \\ 
&  & 
\end{array}%
$$ is a commuting square \ for all $|S_{n}|<\infty $ ,$S_{n}\subset
S_{n+1}\subset I$. 
\end{proof}
\bigskip

\begin{remark}\label{3.7} \mbox{}
\begin{enumerate}
\item
We would like to mention that above Corollary 3.5 recovers Corollary 2.4 in ~\cite{SiSmWhWi} which is the singular version.

\item
 It is worth mentioning some questions that we believe will lead  to a better understanding of these phenomenas:
Is it possible to replace the commuting square condition in the 
Proposition 3.5 by a weaker condition, such that the same conclusion follows? If
yes, what kind of condition? Is it true that we can completely drop the
commuting square condition in the case where $M_{n}$ are hyperfinite factors
for all $n\in \mathbb{N}$?

\item In particular, Proposition 3.5 together with Voiculescu's famous result
~\cite{V1} says: in $L(\mathbb{F}_{n})\overline{\otimes }L(\mathbb{F}_{m})$ we
cannot have a Cartan subalgebra of the form $A\overline{\otimes }B$ with $%
A\subset L(\mathbb{F}_{n})$, $B,\subset L(\mathbb{F}_{m}),$ which was
expected.
\end{enumerate}
\end{remark}

\medskip

 We end this section by presenting a result that estimates the normalizing algebras for certain subalgebras of $II1$ factors arising from cross-product construction.

\def\Aut{{\rm Aut}}

\begin{theorem}\label{3.8}
Let $G$ be a discrete ICC group  and $H$ a subgroup such that the triple $L(H)\subset L(N_{G}(H)\subset L(G)$ satisfies the relative WAHP. Let $N$ be a \II1 factor and let $A\subset N$  be a MASA. 
Suppose there exists an outer action $\alpha :G\rightarrow \Aut(N)$ which satisfies 
that for all $g\in N_{G}(H)$ we have  $\alpha_{g}(A)=A$.

Then,
 for all $g\in N_{G}(H)$ we have 
$\alpha_{g}(\mathcal{N}_{N}(A)^{\prime \prime})=\mathcal{N}_{N}(A)^{\prime \prime}$, and
$$A\rtimes _{\alpha}N_{G}(H)\subseteq \mathcal{N}_{N\rtimes _{\alpha
}G}(A\rtimes _{\alpha }H)^{\prime\prime}\subseteq \mathcal{N}_{N}(A)^{\prime\prime}\rtimes _{\alpha }N_{G}(H).$$
In particular, if $A\subset N$ is a singular MASA, then 
$$A\rtimes _{\alpha
}N_{G}(H)=\mathcal{N}_{N\rtimes _{\alpha }G}(A\rtimes _{\alpha
}H)^{\prime\prime}.$$
\end{theorem}

\begin{proof}

We only prove 
$$(***)\qquad \:\mathcal{N}_{N\rtimes _{\alpha
}G}(A\rtimes _{\alpha }H)^{\prime\prime}\subseteq \mathcal{N}_{N}(A)^{\prime\prime}\rtimes _{\alpha }N_{G}(H),$$
 the other parts being trivial. 
First let us denote by $u_g$ the unitaries that implements the action of $G$ on $N$.

In the same spirit as before, to prove $(***)$ it is enough to check that the triple
$$A\rtimes _{\alpha}H\subseteq \mathcal{N}_{N}(A)^{\prime\prime}\rtimes _{\alpha }N_{G}(H)\subseteq N\rtimes_{\alpha}G$$
satisfies the relative WAHP .

Following the Remark 3.4 this is equivalent to verifying the following: \newline  For every  $S\subset G$ finite subset, every $n_g\in N$ with $g\in S$ and every $ \varepsilon > 0$,   
 there  exists an element $u\in \mathcal{U}(A\rtimes_{\alpha}H)$ such that $$%
(****)\:||E_{A\rtimes _{\alpha }H}(n_{g}u_{g}un_{h}u_{h})-E_{A\rtimes _{\alpha }H}(E_{\mathcal{N}_{N}(A)^{\prime \prime}\rtimes _{\alpha }N_{G}(H)}(n_{g}u_{g})uE_{\mathcal{N}_{N}(A)^{\prime \prime}\rtimes _{\alpha }N_{G}(H)}(n_{h}u_{h}))||_{2}<\varepsilon $$ for
all $g,h\in S.$

First, we fix $ \varepsilon > 0$,  $S\subset G$ finite subset and $n_g\in N$ with $g\in S$ and our goal is to built $u\in \mathcal{U}(A\rtimes_{\alpha}H)$ that will satisfies $(****)$.
By the assumption that 
 the triple $L(H)\subseteq L(N_{G}(H))\subseteq L(G)$ satisfies the relative WAHP, there is a unitary $v \in L(H)$ such that:
$$||E_{L(H)}(u_{g}vu_{h})||_{2}<\frac{\varepsilon }{%
2\underset{g,h\in S }{max}||n_{g}||\mbox{ }||n_{h}||}$$
 for $g$ or $%
h\in S \setminus N_{G}(H).$

If $g\in G\setminus N_{G}(H)$ or $h\in G\setminus N_{G}(H)$ we can evaluate:

$%
||E_{A\rtimes _{\alpha }H}(n_{g}u_{g}avn_{h}u_{h})-E_{A\rtimes _{\alpha }H}(E_{\mathcal{N}_{N}(A)^{\prime \prime}\rtimes _{\alpha }N_{G}(H)}(n_{g}u_{g})avE_{\mathcal{N}_{N}(A)^{\prime \prime}\rtimes _{\alpha }N_{G}(H)}(n_{h}u_{h}))||_{2}$

$=||E_{A\rtimes _{\alpha }H}(n_{g}u_{g}avn_{h}u_{h})||_{2}$

$=||E_{A\rtimes _{\alpha }H}(n_{g}\alpha _{g}(a)u_{g}vu_{h}\alpha
_{h^{-1}}(n_{h}))||_{2}$

$\leq ||E_{N\rtimes _{\alpha }H}(n_{g}\alpha _{g}(a)u_{g}vu_{h}\alpha
_{h^{-1}}(n_{h}))||_{2}$

$=||n_{g}\alpha _{g}(a)E_{N\rtimes _{\alpha }H}(u_{g}vu_{h})\alpha
_{h^{-1}}(n_{h})||_{2}$

$\leq ||n_{g}||\mbox{ }||n_{h}||\mbox{ }||E_{N\rtimes _{\alpha
}H}(u_{g}vu_{h})||_{2}$

$=||n_{g}||\mbox{ }||n_{h}||\mbox{ }||E_{L(H)}(u_{g}vu_{h})||_2$

$<\frac{\varepsilon}{2}$ \ for all $a\in \mathcal{U}(A)$.\newline

We used here that $%
\begin{array}[t]{ccc}
N\rtimes _{\alpha }G & \supset  & L(G) \\ 
\cup  &  & \cup  \\ 
N\rtimes _{\alpha }H & \supset  & L(H) \\ 
&  & 
\end{array}%
$ is a commuting square.\newline

So, when $g\in G\setminus N_{G}(H)$ or $h\in G\setminus N_{G}(H)$ $(****)$ holds true for any unitary of the form $u=av$ with $a\in \mathcal{U}(A)$. \qquad($\delta $)

 When both $g,h\in N_{G}(H)$ we denote by $%
r_{g}=n_{g}-E_{\mathcal{N}_{N}(A)^{\prime \prime}}(n_{g}), r_{h}=n_{h}-E_{\mathcal{N}_{N}(A)^{\prime \prime}}(n_{h})$. Next, we approximate $v$ by a finite sum $\Sigma
_{k\in T}\ v(k)u_{k}$ which satisfies 
$|| v - \Sigma
_{k\in T}\ v(k)u_{k} ||_2 < \frac{\varepsilon}{2}$,   where $T$ is a finite subset of $H$. Since $E_{\mathcal{N}_{N}(A)^{\prime \prime}}(r_{g})=E_{\mathcal{N}_{N}(A)^{\prime \prime}}(r_{h})=0$ and $\mathcal{N}_{N}(A)^{\prime \prime}=\alpha _{g^{-1}}(\mathcal{N}_{N}(A)^{\prime \prime})$,  $%
\mathcal{N}_{N}(A)^{\prime \prime }=\alpha _{k}(\mathcal{N}_{N}(A)^{\prime \prime})$, we obviously get: $$E_{\mathcal{N}_{N}(A)^{\prime \prime}}(\alpha
_{g^{-1}}(r_{g}))=E_{\mathcal{N}_{N}(A)^{\prime \prime}}(\alpha _{k}(r_{h}))=0.$$ By the relative WAHP for the triple $A\subseteq \mathcal{N}_{N}(A)^{\prime \prime}\subseteq N$ there exists $a_{\varepsilon }\in \mathcal{U}(A)$
which satisfies:
$$||E_{A}(\alpha _{g^{-1}}(r_{g})a_{\varepsilon }\alpha
_{k}(r_{h})||_{2}<\frac{\varepsilon }{\#(T)\cdot max_{k\in T}|v(k)|}.$$

At this point we can estimate:

$%
||E_{A\rtimes _{\alpha }H}(n_{g}u_{g}a_{\varepsilon }vn_{h}u_{h})-E_{A\rtimes _{\alpha }H}(E_{\mathcal{N}_{N}(A)^{\prime \prime}\rtimes _{\alpha }N_{G}(H)}(n_{g}u_{g})a_{\varepsilon }vE_{\mathcal{N}_{N}(A)^{\prime \prime}\rtimes _{\alpha }N_{G}(H)}(n_{h}u_{h}))||_{2}$

$=||E_{A\rtimes _{\alpha }H}(r_{g}u_{g}avr_{h}u_{h})||_{2}$

$\leq \frac{\varepsilon}{2}+||\Sigma _{k\in T}$ $v(k)E_{A\rtimes _{\alpha
}H}(r_{g}u_{g}a_{\varepsilon }u_{k}r_{h}u_{h})||_{2}$

$\leq \frac{\varepsilon}{2}+\Sigma _{k\in T}|v(k)|\mbox{ }||E_{A\rtimes _{\alpha
}H}(r_{g}u_{g}a_{\varepsilon }u_{k}r_{h}u_{h})||_{2}$

$=\frac{\varepsilon}{2}+\Sigma _{k\in T}|v(k)|\mbox{ }||E_{A\rtimes _{\alpha
}H}(u_{g}\alpha _{g^{-1}}(r_{g})a_{\varepsilon }\alpha _{k}(r_{h})u_{kh}||_{2}$

$\leq \frac{\varepsilon}{2}+\Sigma _{k\in T}|v(k)|\mbox{ }||E_{A\rtimes _{\alpha
}N_{G}(H)}(u_{g}\alpha _{g^{-1}}(r_{g})a_{\varepsilon }\alpha _{k}(r_{h})u_{kh})||_{2}$

$=\frac{\varepsilon}{2}+\Sigma _{k\in T}|v(k)|\mbox{ }||E_{A}(\alpha
_{g^{-1}}(r_{g})a_{\varepsilon }\alpha _{k}(r_{h})||_{2}.$ \qquad($\rho$) 

Next, we define
\begin{enumerate} 
\item$v_{\varepsilon }=v $ if $S\nsubseteq N_G(H)$ 
\item$v_{\varepsilon }=1$ if  $S\subset N_G(H)$ in which case $T=\{e\}$ and $v(e)=1$ in the computation above.

\end{enumerate}
Finally, the computation ($\rho$) together with ($\delta$) show that $(****)$ holds true for $u=a_{\varepsilon }v_{\varepsilon }\in \mathcal{U}(A\rtimes_{\alpha}H)$, which lead to the desired conclusion.

\end{proof}

\bigskip

To this end, for a better understanding of examples of triples of algebras that satisfiy the relative WAHP we would like to mention a group version for it: 
\begin{proposition} Let $F\leq H\leq G$  groups that satisfy the following:

(P)\quad   For every $S\in G\setminus H$ finite subset there exists $f\in F$ such that 
$gfh\notin F$ for all $g,h\in S$ then the triple $$L(F)\subseteq L(H)\subseteq L(G)$$ satisfies the relative WAHP.
If in addition $F\unlhd H$ is normal, then $\mathcal{N}_{L(G)}{(L(F))}^{\prime \prime}= L(H)$.

\end{proposition}
\begin{remark}\label{3.10} \mbox{}
\begin{enumerate} 
 
\item We believe that both Proposition \ref{3.3} and Proposition \ref{3.8}
follow from a more general statement, but so far we were not able to find the
right setting.
\item Finally, we remark that the relative WAHP gives an easy path to Dixmier results ~\cite{Di} and also recaptures estimates of the normalizing algebra in the situations of free products in ~\cite{Po5,DySiSm,IPePo}.
\end{enumerate}
\end{remark}


\providecommand{\bysame}{\leavevmode\hbox to3em{\hrulefill}\thinspace}
\providecommand{\MR}{\relax\ifhmode\unskip\space\fi MR }
\providecommand{\MRhref}[2]{%
  \href{http://www.ams.org/mathscinet-getitem?mr=#1}{#2}
}
\providecommand{\href}[2]{#2}

\end{document}